# Blow-up of solutions to a $p$-Laplace equation

Yuliya Gorb,[*] Alexei Novikov[†]

November 28, 2011

**Abstract**

Consider two perfectly conducting spheres in a homogeneous medium where the current-electric field relation is the power law. Electric field $E$ blows up in the $L^\infty$-norm as $\delta$, the distance between the conductors tends to zero. We give here a concise rigorous justification of the rate of this blow-up in terms of $\delta$. If the current-electric field relation is linear, similar results were obtained earlier in [5, 20, 29, 30].

## 1 Introduction

Our work is motivated by the issue of material failure initiation. Failure initiation occurs in the zones of high concentrations of extreme electric or current fields, heat fluxes and mechanical loads. Such zones are normally created by external loads amplified by composite microstructure. Therefore, the main focus of this study is on high contrast concentrated composites in which conditions for high field concentration have been created. A $p$-Laplace equation on domains with two spherical inclusions is a prototype set-up for models of two-phase nonlinear composite materials.

We assume that the composite occupies a bounded domain in $\mathbb{R}^d$, $d = 2, 3$. We consider a composite material that consists of a background medium that contains two particles of a different material. The distance $\delta$ between particles is much smaller than their sizes. We also assume that particles are perfectly conducting and the background medium is described by the current-electric field relation

$$\boldsymbol{J} = \sigma|\boldsymbol{E}|^{p-2}\boldsymbol{E}, \quad p > 2,\ p \in \mathbb{N}. \tag{1}$$

Relation (1) describes various physical phenomena which includes deformation theory of plasticity, e.g. [13, 23, 26], where $\boldsymbol{E}$ and $\boldsymbol{J}$ are identified with infinitesimal strain and stress, respectively, nonlinear dielectrics, e.g. [8, 11, 15, 27], where $\boldsymbol{E}$ and $\boldsymbol{J}$ are identified with electric field and current, respectively, and fluid flow, e.g. [1, 24], where $\boldsymbol{E}$ and $\boldsymbol{J}$ are identified with rate of strain and fluid stress, respectively. For definiteness we say that $u$ is the electric potential, and the electric field $\boldsymbol{E} = \nabla u$.

The energy in the thin gaps between neighboring particles of the composite described by $\int \boldsymbol{J} \cdot \boldsymbol{E}$ exhibits singular behavior as $\delta \to 0$ [14], as well electric field $\nabla u$ in the composite. If $p = 2$ then [5] (see also [20, 29, 30]) typically there exists $C > 0$ independent of $\delta$ such that

$$\frac{1}{C\sqrt{\delta}} \leq \|\nabla u\|_{L^\infty} \leq \frac{C}{\sqrt{\delta}} \quad \text{for} \quad d=2,\quad \frac{\log \delta^{-1}}{C} \leq \|\nabla u\|_{L^\infty} \leq C \log \delta^{-1} \quad \text{for} \quad d=3, \tag{2}$$

The main result of this paper is the first asymptotic estimate for any $p > 2$. Namely, we show that typically

$$\lim_{\delta \to 0} \|\nabla u\|_{L^\infty} \delta^\gamma = C,\ C > 0, \quad \text{for} \quad d = 2,\ 3,\ p > 2, \tag{3}$$

where the constant $\gamma = \gamma(p, d) > 0$ and $C$ is explicitly computable. We use the method of sub- and super-solutions to obtain (3). This method is applicable when $p = 2$ as well. In this case our argument is shorter than in [5, 20, 29, 30]. Further it allows to compute the precise limit (3). This is the main contribution of our work. We prove (3) for physically relevant dimensions $d = 2, 3$ and for only two inclusions. By similar methods one may obtain estimates when $d > 3$.

---

[*]Department of Mathematics, University of Houston, Houston, TX, 77204-3008, gorb@math.uh.edu
[†]Department of Mathematics, Pennsylvania State University, University Park, PA, 16802, anovikov@math.psu.edu

Our paper is organized as follows. Next Section 2 provides problem formulation and statements of the main results. Chapter 3 presents proofs of the auxiliary Proposition 2.1. Discussion of the linear case ($p = 2$) and comparison with the previous results are given in Chapter 4. The proof of the main theorem is provided in Chapter 5. In the rest of the introduction we review some earlier results on gradient estimates for high contrast composites.

E.Bonnetier and M.Vogelius [7] studied elliptic regularity for 2-dimensional problems with discontinuous coefficients that model two-phase fiber-reinforced composites with touching fibers. Material properties of its constituents were assumed to be finite and strictly positive. They showed that solutions are $W^{1,\infty}$ for sufficiently smooth boundary data, and conjecture the $\delta^{-1/2}$ blow-up rate (2) for high contrast materials.

Y.Y.Li and M.Vogelius [18] showed $C^{1,\alpha}$-regularity, $0 < \alpha \leq 1/4$ of solutions to elliptic equations that model inhomogeneous materials in $\mathbb{R}^d$. They obtained uniform bounds on $|\nabla u|$ independent of the distance between inhomogeneities. These bounds depended on their sizes, shapes, and material properties. They further conjectured that this gradient exhibits a singular behavior when inhomogeneities are infinitely conducting.

Y.Y.Li and L.Nirenberg [17] extended results in [18] to the vectorial case, namely, to linear systems of elasticity. They obtained $C^{1,\alpha}$ interior estimates on domains in $\mathbb{R}^d$ with non-circular inclusions.

H.Ammari, H.Kang, and M.Lim [4] were the first to investigate the case of close-to-touching regime of particles whose conductivities degenerate, that is, the case of a high contrast composite with perfectly conducting or insulating particles in the background medium of finite conductivity. A lower bound on $|\nabla u|$ was constructed there showing blow up in both perfectly conducting and insulating case. This blow up was proved to be of order $\delta^{-1/2}$ in $\mathbb{R}^2$, where $\delta > 0$ is the distance between two circular particles. In their subsequent work with H.Lee and J.Lee [2] they establish upper and lower bounds on the electric field for close-to-touching regime of two circular particles in $\mathbb{R}^2$ with degenerate conductivities. This study reveals that the blow up of the electric field is of order $\delta^{-1/2}$ and in high contrast regime it occurs at the points of the closest distance between particles. The case of a particle close to boundary is also considered for which similar lower and upper bounds for $|\nabla u|$ are established. Essentially 2-dimensional potential theory techniques are used in [2], and the authors point out the importance of the 3-dimensional case.

In a recent paper [3] H.Ammari et. al. extend results of [2, 4] and decompose the solution into two parts whose gradients are bounded and singular, respectively. This decomposition allows for explicit capturing the gradient's blow up of electric field between circular particles in $\mathbb{R}^2$. Also, a case of nonzero permittivity of inclusions is considered in [3] whose presence is shown to reduce the blow up of the gradient.

K.Yun [29, 30] generalized the blow-up results of [2, 4] to the case of particles of arbitrary sufficiently smooth shape in $\mathbb{R}^2$. Using probabilistic methods it is shown there that the blow up rate of the electric field in composites with arbitrary shape particles is the same as that of disks, that is, of $O(\delta^{-1/2})$. Results [29, 30] were extended by M.Lim and K. Yun [20] to the case of spherical particles in $\mathbb{R}^d$, $d \geq 2$. This is the first result where constants in constructed bounds explicitly contain information about geometry of particles.

E. S. Bao, Y.Y. Li, and B.Yin [5] analyzed a model of a composite with two perfectly conducting particles in a bounded domain and away from the external boundary. The optimal upper and lower bounds for the electric field in a composite when the distance between particles is small were obtained. The results of [5], obtained independently of [20, 29, 30], hold for arbitrary shapes of the particles, that are strictly convex in the neighborhood of the point of the shortest distance between particles, and any dimension $d \geq 2$. E. S. Bao, Y.Y. Li, and B.Yin further generalized their results to the case of $N \geq 2$ perfectly conducting particles and to the case of insulating (zero conductivity) particles in [6].

**Acknowledgements.** The authors thank J. Manfredi, I. Fonseca, G. Leoni, L. Ryzhik and G. Lieberman for helpful discussions on the subject of the paper. Y. G. and A. N. were supported by the NSF grants DMS-1016531 and DMS-0908011, respectively.

## 2  Formulation and Main Results

Let $\Omega \in \mathbb{R}^d$, $d = 2, 3$ be an open bounded domain with $C^{1,\alpha_0}$, $0 < \alpha_0 \leq 1$ boundary $\partial \Omega$. It contains two spherical particles $\mathcal{B}_\delta^1$ and $\mathcal{B}_\delta^2$ centered at $x_i$, $i = 1, 2$ and at the distance $\delta$ from each other (see Figure 1). We

assume
$$\operatorname{dist}(\partial\Omega, \mathcal{B}_\delta^1 \cup \mathcal{B}_\delta^2) \geq K \tag{4}$$
for some $K$ independent of $\delta$. Let $\Omega_\delta$ model the background medium of the composite, that is, $\Omega_\delta = \Omega/\overline{(\mathcal{B}_\delta^1 \cup \mathcal{B}_\delta^2)}$.

Consider
$$\begin{aligned}
\nabla \cdot (|\nabla u_\delta|^{p-2} \nabla u_\delta) &= 0, & \text{in } \Omega_\delta \\
u_\delta &= T_\delta^i & \text{in } \mathcal{B}_\delta^i,\ i=1,2, \\
\int_{\partial \mathcal{B}_\delta^i} |\nabla u_\delta|^{p-2} \boldsymbol{n} \cdot \nabla u_\delta\, ds &= 0, & i=1,2, \\
u_\delta &= U(\boldsymbol{x}), & \text{on } \partial\Omega
\end{aligned} \tag{5}$$

where a bounded weak solution $u_\delta$ represents the electric potential in $\Omega_\delta$, $p \in \mathbb{N}$ and $p \geq 2$, and $U(\boldsymbol{x})$ is the given applied potential on the external boundary $\partial\Omega$. Note $u_\delta$ is constant $T_\delta^i$ on the particle $\mathcal{B}_\delta^i$, $i = 1, 2$, which should be found while solving (5).

In order to formulate our main result for (5) we first describe the meaning of the limit in (3). Given two particles $B^i$ in a domain $\Omega$ consider a family of auxiliary problems where the two particles move along the line connecting their centers until they touch. When $B^i$ are located at the distance $\delta > 0$ from each other we denote them $\mathcal{B}_\delta^i$ and it gives us Figure 1 and (5).

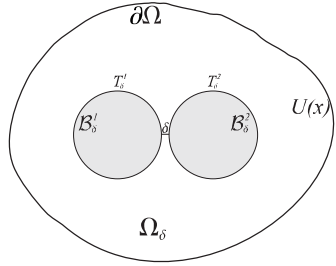

Figure 1: Composite occupying the domain $\Omega$ with particles $\mathcal{B}_\delta^1$ and $\mathcal{B}_\delta^2$

When particles touch at $\delta = 0$ we denote them by $\mathcal{B}_0^i$, see Figure 3. We further construct a *neck* $\Pi_0 = \Pi_0(w)$ of a fixed small width $w > 0$ by cutting out a region that contains the touching point of the two particles, see Fig. 2(a). We denote by $\varsigma_i$, $i = 1, 2$, the part of $\partial \Pi_0$ that lies on the boundary of the corresponding particle. Similarly, one defines a neck $\Pi_\delta = \Pi_\delta(w)$ and $\varsigma_i$, $i = 1, 2$ when particles are $\delta$-distance apart from each other, see Fig. 2(b).

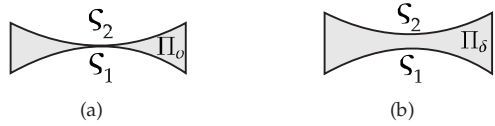

Figure 2: (a) The neck $\Pi_0$; (b) The neck $\Pi_\delta$

Well-posedness of the limiting problem follows, if we look at the variational formulation of (5):
$$u_\delta = \operatorname{argmin}_{u \in U_\delta} \int_{\Omega_\delta} |\nabla u|^p\, d\boldsymbol{x}, \quad U_\delta = \left\{ u \in W^{1,p}(\Omega_\delta):\ u|_{\mathcal{B}_\delta^i} = t_i,\ i=1,2,\ u = U(\boldsymbol{x}) \text{ on } \partial\Omega \right\}. \tag{6}$$

The solution of
$$\begin{aligned}
\nabla \cdot (|\nabla u_0|^{p-2} \nabla u_0) &= 0, & \text{in } \Omega_0, \\
u_0 &= T_0, & \text{in } \mathcal{B}_0^1 \cup \mathcal{B}_0^2, \\
\int_{\partial \mathcal{B}_0^1} |\nabla u_0|^{p-2} \boldsymbol{n} \cdot \nabla u_0\, ds + \int_{\partial \mathcal{B}_0^2} |\nabla u_0|^{p-2} \boldsymbol{n} \cdot \nabla u_0\, ds &= 0, \\
u_0 &= U(\boldsymbol{x}), & \text{on } \partial\Omega.
\end{aligned} \tag{7}$$

is the minimizer of (6) with $\delta = 0$. Note that in this case $u_0$ is the same constant on *both* particles.

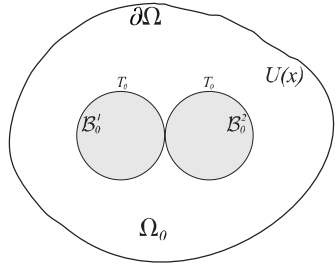

Figure 3: Particles location in the limiting problem (7) whose solution is $u_0$

If $\Omega_\delta$ has a $C^{1,\alpha_0}$ boundary ($0 < \alpha_0 \leq 1$), then [19] there exists a positive constant $\alpha = \alpha(\alpha_0, p, d)$ such that the bounded weak solution of (5) satisfies $u_\delta \in C^{1,\alpha}(\bar{\Omega}_\delta)$. The following proposition characterizes convergence of $u_\delta$ to $u_0$ in $C^{1,\alpha}$.

**Proposition 2.1** *Let $u_\delta$ be the solution of (5), and $u_0$ be the solution of (7). Then there is a constant $\alpha > 0$, $\alpha = \alpha(n, p)$ so that*
$$\lim_{\delta \to 0} \|u_\delta - u_0\|_{C^{1,\alpha}(\mathbf{K})} = 0$$
*for any compact $\mathbf{K} \subset\subset \Omega_0$. Further, for any $i = 1, 2$ and for any neck width $w$*
$$\int_{\partial \mathcal{B}_0^i \setminus \varsigma_i} |\nabla u_0|^{p-2} \boldsymbol{n} \cdot \nabla u_0 = \lim_{\delta \to 0} \int_{\partial \mathcal{B}_\delta^i \setminus \varsigma_i} |\nabla u_\delta|^{p-2} \boldsymbol{n} \cdot \nabla u_\delta. \tag{8}$$

While proving Proposition 2.1 we show that for any neck width $w$ there exists $C = C(w, U)$ so that
$$|\nabla u_\delta(\boldsymbol{x})| \leq C, \text{ for any } \boldsymbol{x} \in \overline{\Omega}_\delta \setminus \overline{\Pi}_\delta. \tag{9}$$

It means that the only possible place for singular behavior of $|\nabla u_\delta|$ is between the particles.

We emphasize that $u_\delta$ and $u_0$ differ in their constraints. Namely, for any $\delta > 0$ the integral of the flux of $u_\delta$ along the boundary of each of the particles is zero, see the third condition in (5). In contrast, for $u_0$ we only assume that the total flux of $u_0$ along the boundary of the both particles is zero, see the third condition in (7). Generically the quantity
$$\mathcal{R}_0 := \int_{\partial \mathcal{B}_0^2} |\nabla u_0|^{p-2} \boldsymbol{n} \cdot \nabla u_0 \, ds \tag{10}$$
is not zero. Since $u_0 \in W^{1,\infty}$, we may use (8) to conclude
$$\left| \mathcal{R}_0 - \lim_{\delta \to 0} \int_{\partial \mathcal{B}_\delta^2 \setminus \varsigma_2} |\nabla u_\delta|^{p-2} \boldsymbol{n} \cdot \nabla u_\delta \, ds \right| \leq Cw,$$
for any any neck width $w$. If $\mathcal{R}_0 \neq 0$ we assume without loss of generality that $\mathcal{R}_0 > 0$. Assumption $\mathcal{R}_0 > 0$ implies that for sufficiently small $\delta$ we have the inequality $T_\delta^2 > T_\delta^1$. It turns out $\mathcal{R}_0 = \mathcal{R}_0[U]$ is the key characteristic parameter of the $W^{1,\infty}$ blow-up of $u_\delta$.

**Theorem 2.2** *Suppose $d \leq p$. If $\mathcal{R}_0 > 0$, then*
$$\lim_{\delta \to 0} (T_\delta^2 - T_\delta^1)^{p-1} \delta^{-\gamma} = \mathcal{R}_0 / C_o, \tag{11}$$
*where $C_o$ is an explicitly computable constant, given in Table (12) below. Similarly, if $d = 3$, $p = 2$, then*
$$\lim_{\delta \to 0} (T_\delta^2 - T_\delta^1)^{p-1} \ln \delta^{-1} = \mathcal{R}_0 / C_o.$$

|         | $d=2$ | $d=3$ |
|---------|-------|-------|
| $p=2$   | $\pi\sqrt{R}$ | $\pi R \ln R$ |
| $p=3$   | $\pi\sqrt{R}/2$ | $\pi R/2$ |
| $p=4$   | $3\pi\sqrt{R}/8$ | $\pi R/8$ |
| general | $\pi\prod_{k=1}^{p-2}\frac{k-1/2}{k}\sqrt{R}$ | $\frac{\pi}{2^{p-1}(p-2)}R$ |

(12)

For the sake of brevity, we focus in our computations on the case $d \leq p$ only, the log-case can be treated similarly.

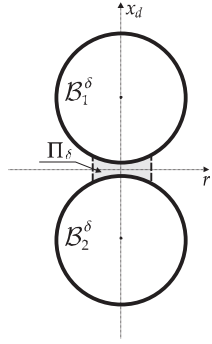

Figure 4: Two disks $\mathcal{B}_\delta^1$ and $\mathcal{B}_\delta^2$, and the local coordinate system. Here $r = x$, $x_d = y$ for $d = 2$, and $r = \sqrt{x^2 + y^2}$, $x_d = z$ for $d = 3$

Suppose we rotate and shift the domain $\Omega_\delta$ so that $\mathcal{B}_\delta^2$ is above $\mathcal{B}_\delta^1$ as depicted on Figure 4. In the proof of Theorem 2.2 we use barriers to show that

$$\frac{T_\delta^2 - T_\delta^1}{\delta + x^2/R}(1 + O(\delta)) - C \leq \boldsymbol{n} \cdot \nabla u_\delta(\boldsymbol{x}) \leq \frac{T_\delta^2 - T_\delta^1}{\delta + x^2/R}(1 + O(\delta)) + C, \quad \boldsymbol{x} \in \varsigma_2, \tag{13}$$

where $C = C(\max(U(x)), \min(U(x)), K)$, $K$ is given in (4). Since $|\nabla u_\delta|$ satisfies the Maximum principle [10], we may use (9), and (13) to immediately obtain the following Corollary of Theorem 2.2.

**Corollary 2.3** *Suppose $d \leq p$. If $u_\delta$ solves (5) then*

$$\|\nabla u_\delta\|_{L^\infty(\Omega_\delta)} = \frac{|T_\delta^2 - T_\delta^1|}{\delta}[1 + O(\delta)] = \left(\frac{\mathcal{R}_0}{C_o}\right)^{\frac{1}{p-1}} \delta^{\frac{\gamma}{p-1} - 1}[1 + O(\delta)]. \tag{14}$$

*Similarly, if $d = 3$, $p = 2$, then*

$$\|\nabla u_\delta\|_{L^\infty(\Omega_\delta)} = \left(\frac{\mathcal{R}_0}{C_o}\right)^{\frac{1}{p-1}} \frac{\left(\ln \delta^{-1}\right)^{\frac{1}{p-1}}}{\delta}[1 + O(\delta)].$$

## 3 Proof of Proposition 2.1

Let $u_\delta$ be the solution of (5). We claim that there exists a constant $C = C(U) > 0$ independent of $\delta$ such that

$$\|\nabla u_\delta\|_{L^\infty(\partial\Omega)} \leq C. \tag{15}$$

We prove (15) by constructing upper and lower barriers for $u$ at $\boldsymbol{x} \in \partial\Omega_0$. The upper barrier $\psi_e$ for $\boldsymbol{x}_0 \in \partial\Omega$ can be chosen as the solution to the following problem:

$$\begin{aligned}
\nabla \cdot (|\nabla\psi_e|^{p-2}\nabla\psi_e) &= 0, & &\text{in } \Omega_\delta \\
\psi_e &= \max_{\partial\Omega} U(x), & &\text{on } \partial\mathcal{B}_\delta^1 \cup \partial\mathcal{B}_\delta^2 \\
\psi_e &= U(x), & &\text{on } \partial\Omega
\end{aligned} \tag{16}$$

It implies there exists such a constant $C_{eu} > 0$ independent of $\delta$ that

$$|\nabla \psi_e(\boldsymbol{x}_0)| \leq C_{eu}. \tag{17}$$

Similarly, a lower barrier $\varphi_e$ for $u$ at $\boldsymbol{x}_0 \in \partial \Omega$ is chosen to be the solution of:

$$\begin{aligned}
\nabla \cdot (|\nabla \varphi_e|^{p-2} \nabla \varphi_e) &= 0, & &\text{in } \Omega_\delta \\
\varphi_e &= \min_{\partial \Omega} U(x)), & &\text{on } \partial \mathcal{B}_\delta^1 \cup \partial \mathcal{B}_\delta^2 \\
\varphi_e &= U(x), & &\text{on } \partial \Omega,
\end{aligned} \tag{18}$$

and there exists such a constant $C_{e\ell} > 0$ independent of $\delta$ that

$$|\nabla \varphi_e(\boldsymbol{x}_0)| \leq C_{e\ell}. \tag{19}$$

Since

$$0 \leq |\boldsymbol{n} \cdot \nabla u_\delta(\boldsymbol{x})| \leq \max\{|\boldsymbol{n} \cdot \nabla \varphi_e(\boldsymbol{x})|, |\boldsymbol{n} \cdot \nabla \psi_e(\boldsymbol{x})|\} \quad \text{if } \boldsymbol{x}_0 \in \partial \Omega,$$

estimates (17), (19) imply (15).

We now show that for any fixed neck width $w$ there exists a constant $C = C(w)$ such that

$$|\nabla u_\delta(x)| \leq C, \quad \text{for} \quad x \in \partial \Pi_\delta^\pm. \tag{20}$$

Consider an auxiliary function

$$\Phi(\boldsymbol{x}) = f(x)|\nabla u_\delta|^2 + \kappa u_\delta^2, \quad \boldsymbol{x} = (x, y) \in \Omega_\delta.$$

This proof focuses on 2D, the 3D case follows similarly. Choose function $f(x)$ defined on the coordinate system as on Fig. 4 to be as follows:

$$f(x) = \begin{cases} 0, & |x| \leq \dfrac{w}{2}, \\ \dfrac{8}{w^2}\left(|x| - \dfrac{w}{2}\right)^2, & \dfrac{w}{2} < |x| < \dfrac{3}{4}w, \\ -\dfrac{8}{w^2}(|x| - w)^2 + 1, & \dfrac{3}{4}w^2 < |x| < w, \\ 1, & |x| \geq w. \end{cases} \tag{21}$$

There exists a constant $\kappa = O(1/w^2)$ so that $\Phi$ satisfies the maximum principle [9, 22]. Then

$$\max_{\partial \Pi_\delta^\pm} |\nabla u_\delta|^2 \leq \max_{\Omega_\delta} \Phi(x, y) \leq C(w) + \max_{\partial \mathcal{B}_\delta^i, \partial \Omega} f(x)|\nabla u_\delta|^2. \tag{22}$$

Since $|\nabla u_\delta| = O(1)$ on $\partial \Omega$ by (15), we only need to obtain estimates for $f(x)|\nabla u_\delta|^2$ on $\partial \mathcal{B}_\delta^i$, $i = 1, 2$ for $\frac{1}{2}w < |x|$. This can be shown by constructing the upper and the lower barrier for $u_\delta$ at $\boldsymbol{x} \in \varsigma_i$, $i = 1, 2$. It is straightforward to verify that

$$\psi(\boldsymbol{x}) = \begin{cases} a|\boldsymbol{x} - \boldsymbol{x}_0|^\beta + b, & d \neq p \\ a \log|\boldsymbol{x} - \boldsymbol{x}_0| + b, & d = p \end{cases} \tag{23}$$

is the solution of the $p$-Laplace equation

$$\nabla \cdot (|\nabla \psi|^{p-2} \nabla \psi) = 0, \text{ on } \mathbb{R}^d \setminus \boldsymbol{x}_0,$$

for any $\boldsymbol{x}_0$, $a$, $b$, and $\beta = (p-d)/(p-1)$. Once again, for brevity we focus in our computations on the case $d \neq p$ only, the log-case can be treated similarly.

The upper barrier for $u_\delta$ at $\boldsymbol{x} = (x, y) \in \varsigma_1$ is obtained as follows. Construct a circle $\mathcal{C}_1$ of the radius $r_1 = w/100$ (see Fig. 7) centered at a point inside $\mathcal{B}_\delta^1$ that touches $\partial \mathcal{B}_\delta^1$ at $\boldsymbol{x} \in \varsigma_1$. Construct a concentric to $\mathcal{C}_1$ circle denoted by $\mathcal{C}_2$ of a radius $r_2 = r_2(r_1)$ that touches $\partial \mathcal{B}_\delta^2$, and

$$r_2 = \delta + r_1 + \frac{1}{2R}\left(\frac{R - r_1}{R}\right)\left(\frac{2R - r_1}{R}\right)x^2. \tag{24}$$

Then an upper barrier $\overline{v}$ for $u_\delta$ at $\boldsymbol{x} \in \varsigma_1$ is obtained using the radial solution (23) that satisfies: $\overline{v}(r_1) = T_\delta^1$, $\overline{v}(r_2) = \max_{\partial\Omega} U(x)$, and the lower barrier $\underline{v}$ for for $u_\delta$ at $\boldsymbol{x} \in \varsigma_1$ is obtained using the radial solution (23) that satisfies: $\underline{v}(r_1) = T_\delta^1$, $\underline{v}(r_2) = \min_{\partial\Omega} U(x)$. Then

$$|\nabla u_\delta(\boldsymbol{x})| = |(\nabla u_\delta \cdot \boldsymbol{n})| \leq \max\left\{C, |\nabla \overline{v}|, |\nabla \underline{v}|\right\}, \quad \boldsymbol{x} \in \varsigma_1, \tag{25}$$

where

$$|\nabla \overline{v}| = \left|\frac{d\overline{v}}{dr}\right| = |\max_{\partial\Omega} U(x) - T_\delta^1| \frac{\beta r_1^{\beta-1}}{r_2^\beta - r_1^\beta}, \quad |\nabla \underline{v}| = \left|\frac{d\underline{v}}{dr}\right| = |T_\delta^1 - \min_{\partial\Omega} U(x)| \frac{\beta r_1^{\beta-1}}{r_2^\beta - r_1^\beta}.$$

Since $r_2 = r_1 + O(w)$ then $r_2^\beta - r_1^\beta = O(w)$, hence, we have

$$|\nabla u_\delta(\boldsymbol{x})| \leq C, \quad C = C(w), \quad \boldsymbol{x} = (x,y) \in \varsigma_1, \quad \text{and} \quad \frac{w}{2} < |x| < R.$$

Similarly, one can bound the gradient at $\boldsymbol{x} \in \varsigma_2$. Substituting these estimates back to (22) we have that there exists $C = C(w)$ such that (20) holds.

The maximum principle and (22) imply (9) that states

$$|\nabla u_\delta(\boldsymbol{x})| \leq C, \text{ for any } \boldsymbol{x} \in \overline{\Omega}_\delta \setminus \overline{\Pi}_\delta. \tag{26}$$

Now the estimate (26) means that, up to a subsequence, $u_\delta \to u_*$, as $\delta \to 0$ strongly in $W^{1,q}$ for any $q < \infty$. We claim that $u_0 = u_*$, where $u_0$ solves (7). Indeed the only condition that needs to be verified is:

$$\int_{\partial \mathcal{B}_0^1} |\nabla u_*|^{p-2} \boldsymbol{n} \cdot \nabla u_* \, ds + \int_{\partial \mathcal{B}_0^2} |\nabla u_*|^{p-2} \boldsymbol{n} \cdot \nabla u_* \, ds = 0. \tag{27}$$

For any small parameter $\varepsilon > 0$ construct a domain $\mathcal{K}_\varepsilon$ with $C^\infty$ boundary such that it approximates $\mathcal{B}_0^1 \cup \mathcal{B}_0^2$ arbitrary well, that is, each point on $\partial \mathcal{K}_\varepsilon$ is located at a distance $O(\varepsilon)$ from $\partial(\mathcal{B}_0^1 \cup \mathcal{B}_0^2)$, see Figure 5.

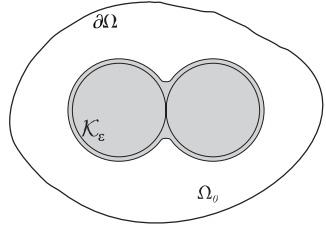

Figure 5: Smoothening $\mathcal{K}_\varepsilon$ of the domain $\mathcal{B}_0^1 \cup \mathcal{B}_0^2$

From [28] we know that there is a constant $\alpha > 0$, $\alpha = \alpha(n, p, \alpha_0)$ so that

$$\|u_\delta\|_{C^{1,\alpha}} \leq C$$

on $\overline{\Omega_0/\mathcal{K}_\varepsilon}$ uniformly in sufficiently small $\delta$. We, therefore, obtain that $u_\delta \to u_*$ strongly in $C^1(\overline{\Omega_0/\mathcal{K}_\varepsilon})$. Integrating the first equation in (5) over $\mathcal{K}_\varepsilon/(\mathcal{B}_0^1 \cup \mathcal{B}_0^2)$, and using the third equation in (5) we obtain

$$\int_{\partial \mathcal{K}_\varepsilon} |\nabla u_\delta|^{p-2} \boldsymbol{n} \cdot \nabla u_\delta \, ds = 0. \tag{28}$$

Since $u_\delta \to u_*$ strongly in $C^1(\overline{\Omega_0/\mathcal{K}_\varepsilon})$ we have

$$\int_{\partial \mathcal{K}_\varepsilon} |\nabla u_*|^{p-2} \boldsymbol{n} \cdot \nabla u_* \, ds = 0.$$

Since $\varepsilon$ is arbitrary we obtain (27). Using [28] $u_0$ is the unique $C^{1,\alpha}(\Omega_0)$ solution of (7). Uniqueness of $u_0$ allows us to conclude that there is a point-wise convergence $u_\delta \to u_* = u_0$ for all $\delta \to 0$. From [28] we also conclude that the convergence is in $C^{1,\alpha}(\mathbf{K})$, $\mathbf{K} \subset\subset \Omega$ as $\delta \to 0$.

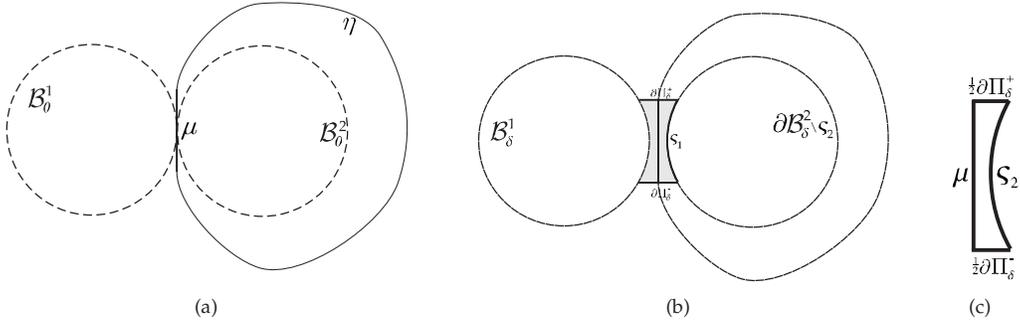

Figure 6: (a) The smooth curve $\Gamma = \mu \cup \eta$ and its flat piece $\mu$; (b) The portion $\varsigma_2$ of $\partial \mathcal{B}_\delta^2$ that belongs to the neck $\Pi_\delta$; (c) The half-neck $\frac{1}{2}\Pi_\delta$

Proof of (8) follows, once we recall that outside any neck $u_\delta$ and $u_0$ are uniformly bounded in $W^{1,\infty}$. So we can take a family of contours depicted on Figure 6; integrate the $p$-Laplacian in the interior of these contours; conclude that the integrals over the boundary of the particles are equal (up to a small constant) to integrals in the interior, where we have convergence. Let us provide more detail.

We first construct a closed $C^\infty$-curve $\Gamma$ that contains $\mathcal{B}_0^2$ and passes through the point where two circles $\mathcal{B}_0^2$ and $\mathcal{B}_0^1$ touch. This curve is large enough to contain $\mathcal{B}_\delta^2$ inside when particles are "moved away" at the distance $\delta$ between each other (recall, in such a configuration particles are denoted by $\mathcal{B}_\delta^1$ and $\mathcal{B}_\delta^2$). We also assume that this curve has a flat portion of the length $2w$ denoted by $\mu$ which is perpendicular to the line connecting the two particles, and the rest of the curve is $\eta = \Gamma \setminus \mu$, see Fig. 6(a). Then

$$\mathcal{R}_0 = \int_\mu |\nabla u_0|^{p-2} (\boldsymbol{n} \cdot \nabla u_0)\, ds + \int_\eta |\nabla u_0|^{p-2} (\boldsymbol{n} \cdot \nabla u_0)\, ds.$$

Since $|\nabla u_0| < C$ everywhere in $\Omega_0$ [7] we have

$$\left| \mathcal{R}_0 - \int_\eta |\nabla u_0|^{p-2} (\boldsymbol{n} \cdot \nabla u_0)\, ds \right| = \left| \int_\mu |\nabla u_0|^{p-2} (\boldsymbol{n} \cdot \nabla u_0)\, ds \right| \leq Cw. \tag{29}$$

We now "move away" particles on the distance $\delta$ between each other, see Fig. 6(b), and consider a region between the two curves: $\partial \mathcal{B}_\delta^2$ and $\Gamma$ where we still have $\nabla \cdot |\nabla u_\delta|^{p-2} \nabla u_\delta = 0$. Multiplying this equation by $u_\delta$ and integrating by parts we obtain

$$\int_\mu |\nabla u_\delta|^{p-2} (\boldsymbol{n} \cdot \nabla u_\delta)\, ds + \int_\eta |\nabla u_\delta|^{p-2} (\boldsymbol{n} \cdot \nabla u_\delta)\, ds - \int_{\partial \mathcal{B}_\delta^2} |\nabla u_\delta|^{p-2} (\boldsymbol{n} \cdot \nabla u_\delta)\, ds = 0, \tag{30}$$

where the third integral of the left hand side is zero. This integral we split into the piece that belongs to the neck $\Pi_\delta$ denoted as above by $\varsigma_2$ and the rest of the boundary $\partial \mathcal{B}_\delta^2 \setminus \varsigma_2$. Then we have

$$\int_{\partial \mathcal{B}_\delta^2 \setminus \varsigma_2} |\nabla u_\delta|^{p-2} (\boldsymbol{n} \cdot \nabla u_\delta)\, ds = \int_\mu |\nabla u_\delta|^{p-2} (\boldsymbol{n} \cdot \nabla u_\delta)\, ds + \int_\eta |\nabla u_\delta|^{p-2} (\boldsymbol{n} \cdot \nabla u_\delta)\, ds - \int_{\varsigma_2} |\nabla u_\delta|^{p-2} (\boldsymbol{n} \cdot \nabla u_\delta)\, ds. \tag{31}$$

Since we have a pointwise convergence of $u_\delta$ to $u_0$ in any compact $K \subset\subset \Omega \setminus (\mathcal{B}_0^1 \cup \mathcal{B}_0^2)$ due to [28] we have that

$$\int_\eta |\nabla u_\delta|^{p-2} (\boldsymbol{n} \cdot \nabla u_\delta)\, ds \to \int_\eta |\nabla u_0|^{p-2} (\boldsymbol{n} \cdot \nabla u_0)\, ds, \quad \text{as} \quad \delta \to 0. \tag{32}$$

Now we multiply the equation $\nabla \cdot |\nabla u_\delta|^{p-2} \nabla u_\delta = 0$ in the region between $\mu$ and $\varsigma_2$ and integrate by parts. We obtain that

$$\int_\mu |\nabla u_\delta|^{p-2} (\boldsymbol{n} \cdot \nabla u_\delta)\, ds - \int_{\varsigma_2} |\nabla u_\delta|^{p-2} (\boldsymbol{n} \cdot \nabla u_\delta)\, ds + \int_{\frac{1}{2}\partial \Pi_\delta^\pm} |\nabla u_\delta|^{p-2} (\boldsymbol{n} \cdot \nabla u_\delta)\, ds = 0, \tag{33}$$

where $\frac{1}{2}\partial\Pi_\delta^\pm$, see Fig. 6(c), is the half of the lateral boundary of the neck attached to $\mathcal{B}_\delta^2$. So, now substituting (33) into (31) and using (32) we have

$$\left| \int_{\partial\mathcal{B}_\delta^2\setminus\varsigma_2} |\nabla u_\delta|^{p-2}(\boldsymbol{n}\cdot\nabla u_\delta)\,ds - \int_{\partial\mathcal{B}_0^2\setminus\varsigma_2} |\nabla u_0|^{p-2}\boldsymbol{n}\cdot\nabla u_0\,ds \right|$$

$$= \left| \int_\eta |\nabla u_\delta|^{p-2}(\boldsymbol{n}\cdot\nabla u_\delta)\,ds - \int_{\frac{1}{2}\partial\Pi_\delta^\pm} |\nabla u_\delta|^{p-2}(\boldsymbol{n}\cdot\nabla u_\delta)\,ds - \int_{\partial\mathcal{B}_0^2\setminus\varsigma_2} |\nabla u_0|^{p-2}\boldsymbol{n}\cdot\nabla u_0\,ds \right|$$

$$\leq \left| \int_\eta |\nabla u_0|^{p-2}(\boldsymbol{n}\cdot\nabla u_0)\,ds - \int_{\frac{1}{2}\partial\Pi_\delta^\pm} |\nabla u_\delta|^{p-2}(\boldsymbol{n}\cdot\nabla u_\delta)\,ds - \int_{\partial\mathcal{B}_0^2\setminus\varsigma_2} |\nabla u_0|^{p-2}\boldsymbol{n}\cdot\nabla u_0\,ds \right| + o(\delta)$$

$$\leq \left| \int_{\frac{1}{2}\partial\Pi_\delta^\pm} |\nabla u_\delta|^{p-2}(\boldsymbol{n}\cdot\nabla u_\delta)\,ds - \int_{\frac{1}{2}\partial\Pi_\delta^\pm} |\nabla u_0|^{p-2}(\boldsymbol{n}\cdot\nabla u_0)\,ds \right| + o(\delta) \to 0, \quad \text{as} \quad \delta \to 0.$$

This estimate leads to (8). Finally, substituting (32), (33) into (31) and using (29) we obtain

$$\lim_{\delta\to 0}\left| \mathcal{R}_0 - \int_{\partial\mathcal{B}_\delta^2\setminus\varsigma_2} |\nabla u_\delta|^{p-2}(\boldsymbol{n}\cdot\nabla u_\delta)\,ds \right| = \left| \mathcal{R}_0 - \int_{\partial\mathcal{B}_0^2\setminus\varsigma_2} |\nabla u_0|^{p-2}(\boldsymbol{n}\cdot\nabla u_0)\,ds \right|$$

$$= \left| \int_{\varsigma_2} |\nabla u_0|^{p-2}(\boldsymbol{n}\cdot\nabla u_0)\,ds \right| \leq Cw.$$

□

## 4  Linear case $p = 2$

For $p = 2$ we chose to compare our results with [5], because it is easy to interpret in terms of our problem. As it is mentioned in [5] the results [2, 4, 20, 29, 30] are essentially the same. Consider a functional

$$Q_\delta[U] = \int_{\partial\mathcal{B}_\delta^1} \frac{\partial v_3}{\partial n} \int_{\partial\Omega} \frac{\partial v_2}{\partial n} - \int_{\partial\mathcal{B}_\delta^2} \frac{\partial v_3}{\partial n} \int_{\partial\Omega} \frac{\partial v_1}{\partial n}, \qquad a_{ij} = \int_{\partial\mathcal{B}_\delta^i} \frac{\partial v_j}{\partial n}$$

where $v_i$ ($i = 1, 2, 3$) solve the following problems:

$$\begin{aligned}
\triangle v_1 &= 0, \quad \text{in } \Omega_\delta & \triangle v_2 &= 0, \quad \text{in } \Omega_\delta & \triangle v_3 &= 0, \quad \text{in } \Omega_\delta \\
v_1 &= 1, \quad \text{on } \partial\mathcal{B}_\delta^1 & v_2 &= 1, \quad \text{on } \partial\mathcal{B}_\delta^2 & v_3 &= 0, \quad \text{on } \partial\mathcal{B}_\delta^1 \cup \partial\mathcal{B}_\delta^2 \\
v_1 &= 0, \quad \text{on } \partial\mathcal{B}_\delta^2 \cup \partial\Omega & v_2 &= 0, \quad \text{on } \partial\mathcal{B}_\delta^1 \cup \partial\Omega & v_3 &= U(x), \quad \text{on } \partial\Omega
\end{aligned}$$

The main result, Theorems 1.1 and 1.2 of [5] states that

$$C|Q_\delta[U]|\sqrt{\delta} \leq |T_\delta^2 - T_\delta^1| \leq C(U)\sqrt{\delta}, \quad C\frac{|Q_\delta[U]|}{|\ln\delta|} \leq |T_\delta^2 - T_\delta^1| \leq \frac{C(U)}{|\ln\delta|}$$

for $d = 2, 3$, respectively. Thus $Q_\delta[U]$ is a characteristic parameter, that determines blow-up. Let us relate it to our $\mathcal{R}_0$. We verify, that

$$Q_\delta[U] = \int_{\partial\mathcal{B}_\delta^1} \frac{\partial v_\delta}{\partial n} \int_{\partial\Omega} \frac{\partial v_2}{\partial n} - \int_{\partial\mathcal{B}_\delta^2} \frac{\partial v_\delta}{\partial n} \int_{\partial\Omega} \frac{\partial v_1}{\partial n},$$

where $v_i$ ($i = 1, 2$) are as above and $v_\delta$ solves

$$\begin{aligned}
\nabla\cdot|\nabla v_\delta|^{p-2}\nabla v_\delta &= 0, & \text{in } \Omega_\delta \\
v_\delta &= T_\delta, & \text{in } \mathcal{B}_\delta^1 \cup \mathcal{B}_\delta^2 \\
\int_{\partial\mathcal{B}_\delta^1} |\nabla v_\delta|^{p-2}\boldsymbol{n}\cdot\nabla v_\delta\,ds + \int_{\partial\mathcal{B}_\delta^2} |\nabla v_\delta|^{p-2}\boldsymbol{n}\cdot\nabla v_\delta\,ds &= 0, & \\
v_\delta &= U(\boldsymbol{x}), & \text{on } \partial\Omega
\end{aligned} \qquad (34)$$

with $p = 2$. Indeed, we have that
$$v_\delta = T_\delta [v_1 + v_2] + v_3 \quad \text{in } \Omega_\delta,$$
then
$$\begin{aligned}
&\int_{\partial \mathcal{B}_\delta^1} \frac{\partial v_\delta}{\partial n} \int_{\partial \Omega} \frac{\partial v_2}{\partial n} - \int_{\partial \mathcal{B}_\delta^2} \frac{\partial v_\delta}{\partial n} \int_{\partial \Omega} \frac{\partial v_1}{\partial n} \\
&= T_\delta \left[ \int_{\partial \mathcal{B}_\delta^1} \frac{\partial v_1}{\partial n} + \int_{\partial \mathcal{B}_\delta^1} \frac{\partial v_2}{\partial n} \right] \int_{\partial \Omega} \frac{\partial v_2}{\partial n} - T_\delta \left[ \int_{\partial \mathcal{B}_\delta^2} \frac{\partial v_1}{\partial n} + \int_{\partial \mathcal{B}_\delta^2} \frac{\partial v_2}{\partial n} \right] \int_{\partial \Omega} \frac{\partial v_1}{\partial n} + Q_\delta[U] \\
&= T_\delta(a_{11} + a_{12})(-a_{12} - a_{22}) - T_\delta(a_{21} + a_{22})(-a_{11} - a_{21}) + Q_\delta[U] = Q_\delta[U],
\end{aligned}$$

using $a_{12} = a_{21}$. Since
$$\int_{\partial \mathcal{B}_\delta^1} \frac{\partial v_\delta}{\partial n} = -\int_{\partial \mathcal{B}_\delta^2} \frac{\partial v_\delta}{\partial n}$$

we conclude
$$Q_\delta[U] = -\left( \int_{\partial \Omega} \frac{\partial v_2}{\partial n} + \int_{\partial \Omega} \frac{\partial v_1}{\partial n} \right) \mathcal{R}_\delta[U]$$

where
$$\mathcal{R}_\delta[U] := \int_{\partial \mathcal{B}_\delta^2} |\nabla v_\delta|^{p-2} \bm{n} \cdot \nabla v_\delta \, ds, \quad p = 2. \tag{35}$$

By the maximum principle
$$0 < C_1 \le -\left( \int_{\partial \Omega} \frac{\partial v_2}{\partial n} + \int_{\partial \Omega} \frac{\partial v_1}{\partial n} \right) \le C_2,$$

therefore asymptotic behavior of $Q_\delta[U]$ is the same as asymptotic behavior of $\mathcal{R}_\delta[U]$. Using methods, as in Section 2.1, we can show the following.

**Proposition 4.1** *Let $v_\delta$ be the solution of (34), and $u_0$ be the solution of (7). Then there is a constant $\alpha > 0$, $\alpha = \alpha(d, p)$ so that*
$$\lim_{\delta \to 0} ||v_\delta - u_0||_{C^{1,\alpha}(\mathbf{K})} = 0$$

*for any compact $\mathbf{K} \subset\subset \Omega_0$. Further*
$$\lim_{\delta \to 0} \mathcal{R}_\delta = \mathcal{R}_0. \tag{36}$$

# 5 Proof of Theorem 2.2

We will use here the method of barriers and the radial solution (23) to prove an estimate, which is slightly stronger than (13).

**Lemma 5.1** *For any constants $T_\delta^1$ and $T_\delta^2$, $m \le T_\delta^1 \le T_\delta^2 \le M$ there exists $\delta_0$ such that for any $\delta \le \delta_0$ the solution of*
$$\begin{aligned}
\nabla \cdot (|\nabla u_\delta|^{p-2} \nabla u_\delta) &= 0, & \text{in } \Omega_\delta \\
u_\delta &= T_\delta^i & \text{in } \mathcal{B}_\delta^i, \ i = 1, 2, \\
u_\delta &= U(\bm{x}), & \text{on } \partial \Omega
\end{aligned} \tag{37}$$

*satisfies an estimate*
$$\frac{T_\delta^2 - T_\delta^1}{\delta + x^2/R}(1 + O(\delta)) - C \le \bm{n} \cdot \nabla u_\delta(\bm{x}) \le \frac{T_\delta^2 - T_\delta^1}{\delta + x^2/R}(1 + O(\delta)) + C, \quad \bm{x} \in \varsigma_2, \tag{38}$$

*where $C = C(\max(U(x)), \min(U(x)), K)$, $K$ given in (4).*

In contrast to (5), the constants $T_\delta^1$ and $T_\delta^2$ in (37) are arbitrary. It implies the solution of (37) may not satisfy the integral identities the flux of $u_\delta$ on $\partial \mathcal{B}_\delta^i$ as in (5).

**Proof of Lemma.**

An upper barrier $\overline{v}$ for $u_\delta$ at the point $(x,y) \in \varsigma_1$ is constructed as follows. Consider a circle of radius $0 < r_1 < R$ that touches the circle $\partial \mathcal{B}_\delta^1$ from within (as in Fig.7) at the point $(x,y)$. Also, another circle of radius $r_2 > r_1$ is considered that centered at the same point as the other one and touches the boundary $\partial \mathcal{B}_\delta^2$. Then

$$r_2 = r_2(x, r_1) = \delta + r_1 + \frac{1}{2}\left(1 - \frac{r_1}{R}\right)\left(2 - \frac{r_1}{R}\right)\frac{x^2}{R} \tag{39}$$

Choosing $r_1 = \delta$ we obtain

$$r_2 - r_1 = \left(\delta + \frac{x^2}{R}\right)(1 + O(\delta))$$

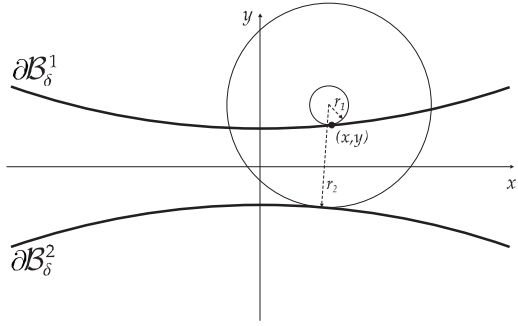

Figure 7: Upper barrier construction at the point $(x,y) \in \varsigma_1$

Then we construct an upper barrier $\overline{v}$ for $u_\delta$ at the point $(x,y) \in \varsigma_1$ via radial solution (23) that takes values

$$\overline{v}(r_1) = T_\delta^1, \quad \overline{v}(r_2) = T_\delta^2.$$

Then

$$\boldsymbol{n} \cdot \nabla u_\delta(x,y) \leq \frac{d\overline{v}}{dr}(r_1) = \beta r_1^{\beta-1}\frac{T_\delta^2 - T_\delta^1}{r_2^\beta - r_1^\beta}. \tag{40}$$

By the mean-value property

$$r_2^\beta - r_1^\beta = \beta r_0^{\beta-1}(r_2 - r_1), \; r_1 < r_0 < r_2$$

and therefore

$$\boldsymbol{n} \cdot \nabla u_\delta(x,y) \leq \beta r_1^{\beta-1}\frac{T_\delta^2 - T_\delta^1}{r_2^\beta - r_1^\beta} = \left(\frac{r_1}{r_0}\right)^{\beta-1}\frac{T_\delta^2 - T_\delta^1}{r_2 - r_1} \leq \frac{T_\delta^2 - T_\delta^1}{r_2 - r_1} = \frac{T_\delta^2 - T_\delta^1}{\delta + x^2/R}(1 + O(\delta)).$$

We can find $\delta_0$ and $C = C(K, m, M, w, \delta_0)$ so that if $(T_\delta^2 - T_\delta^1)/(r_2 - r_1) \geq s$ and $\delta \leq \delta_0$ then

$$u_\delta(\boldsymbol{x}) \leq \overline{v}(\boldsymbol{x}) \text{ for all } \boldsymbol{x} \in \partial\Omega,$$

Thus $\overline{v}$ is an upper barrier if $(T_\delta^2 - T_\delta^1)/(r_2 - r_1) \geq s$ and $\delta \leq \delta_0$. The upper bound in (38) now follows since for other values of $(T_\delta^2 - T_\delta^1)/(r_2 - r_1)$ and $\delta$ we estimate

$$\boldsymbol{n} \cdot \nabla u_\delta(x,y) \leq C.$$

For the lower barrier $\underline{v}$ at $(x,y) \in \varsigma_1$ we consider a circle of a small radius $\rho_1$ that touches the circle $\partial \mathcal{B}_\delta^2$ from within and whose center is located on the line connecting the point $(x,y)$ and the center of $\mathcal{B}_\delta^1$ (as in Fig.8). Denote the distance from the center $(\xi, -\eta)$ of this constructed circle of radius $\rho_1$ to the point $(x,y)$ by $\rho_2$. Then

$$\xi = \frac{R + \rho_2}{R}x, \quad (R + \rho_2)^2 = \xi^2 + (R + \frac{\delta}{2} + \eta)^2, \quad (R - \rho_1)^2 = \xi^2 + (R + \frac{\delta}{2} - \eta)^2$$

we have

$$\rho_2 = \rho_2(x, \rho_1) = -R+(2R+\delta)\left[\sqrt{1-\frac{x^2}{R^2}} - \sqrt{\left(\frac{R-\rho_1}{2R+\delta}\right)^2 - \frac{x^2}{R^2}}\right] = \delta+\rho_1+\frac{2R+\delta}{2R}\left[\frac{2R+\delta}{R-\rho_1}-1\right]\frac{x^2}{R}+O\left(\frac{x^4}{R^2}\right).$$

Choosing $\rho_1 = \delta$ we obtain

$$\rho_2 - \rho_1 \leq \begin{cases} \left(\delta + \frac{x^2}{R}\right)(1+O(\delta)), & |x| \leq \delta^{1/4}, \\ C, & |x| \geq \delta^{1/4}. \end{cases} \tag{41}$$

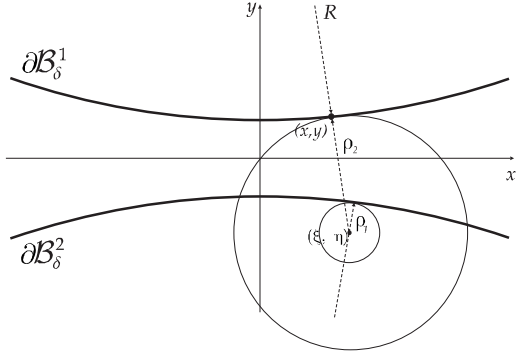

Figure 8: Lower barrier construction at the point $(x,y) \in \varsigma_1$

Then we construct a lower barrier $\underline{v}$ for $u_\delta$ at the point $(x,y) \in \varsigma_1$ via radial solution (23) that takes values

$$\underline{v}(\rho_1) = T_\delta^2, \quad \underline{v}(\rho_2) = T_\delta^1.$$

Then

$$\boldsymbol{n} \cdot \nabla u_\delta(x,y) \geq -\frac{d\underline{v}}{dr}(\rho_2) = \beta\rho_2^{\beta-1}\frac{T_\delta^2 - T_\delta^1}{\rho_2^\beta - \rho_1^\beta} \geq \frac{T_\delta^2 - T_\delta^1}{\rho_2 - \rho_1}. \tag{42}$$

Using (41) in (42) we obtain

$$\boldsymbol{n} \cdot \nabla u_\delta(x,y) \geq -\frac{d\underline{v}}{dr}(\rho_2) \geq \max\left(\frac{T_\delta^2 - T_\delta^1}{\delta + x^2/R}(1+O(\delta)), C\right).$$

As with the upper bound we can find $\delta_0$ and $C = C(K, m, M, w, \delta_0)$ so that if $(T_\delta^2 - T_\delta^1)/(\rho_2 - \rho_1) \geq s$ and $\delta \leq \delta_0$ then

$$u_\delta(\boldsymbol{x}) \geq \underline{v}(\boldsymbol{x}) \text{ for all } \boldsymbol{x} \in \partial\Omega,$$

Thus $\underline{v}$ is an lower barrier if $(T_\delta^2 - T_\delta^1)/(\rho_2 - \rho_1) \geq s$ and $\delta \leq \delta_0$. The lower bound in (38) now follows since for other values of $(T_\delta^2 - T_\delta^1)/(\rho_2 - \rho_1)$ and $\delta$ we again estimate

$$\boldsymbol{n} \cdot \nabla u_\delta(x,y) \geq -C.$$

□

We now prove Theorem 2.2. From (38) we estimate

$$(1+O(\delta))\int_{-w}^{w}\left(\frac{T_\delta^2 - T_\delta^1}{\delta + x^2/R}\right)^{p-1}dx - Cw \leq \int_{\varsigma_1}|\nabla u_\delta|^{p-2}\boldsymbol{n} \cdot \nabla u_\delta \, ds \leq (1+O(\delta))\int_{-w}^{w}\left(\frac{T_\delta^2 - T_\delta^1}{\delta + x^2/R}\right)^{p-1}dx + Cw. \tag{43}$$

Since [22] for any $w > 0$

$$\lim_{\delta \to 0} \delta^\gamma \int_{-w}^{w}\frac{dx}{(\delta + x^2/R)^{p-1}} = C_o$$

the claim of Theorem 2.2 follows. □

# References

[1] Antontsev, S. N., Rodrigues, J. F. : On stationary thermo-rheological viscous flows. *Ann. Univ. Ferrara Sez. VII Sci. Mat.*, **52(1)**, 2006, pp. 19-36.

[2] Ammari, H., Kang, H., Lee, H., Lee, J., Lim, M. : Optimal Bounds on the Gradient of Solutions to Conductivity Problems, *J. Math. Pures Appl.*, **88**, 2007, pp. 307–324.

[3] Ammari, H., Kang, H., Lee, H., Lim, M., Zribi, H. : Decomposition theorems and fine estimates for electrical fields in the presence of closely located circular inclusions, *J. Diff. Eqs.*, **247**, 2009, pp. 2897–2912.

[4] Ammari, H., Kang, H., Lim, M. : Gradient Estimates for Solutions to the Conductivity Problem, *Math. Ann.*, **332:2**, 2005, pp. 277–286.

[5] Bao, E. S., Li, Y. Y., Yin, B. : Gradient Estimates for the Perfect Conductivity Problem, *Arch. Rat. Mech. Anal.*, **193**, 2009, pp. 195–226.

[6] Bao, E. S., Li, Y. Y., Yin, B. : Gradient Estimates for the Perfect and Insulated Conductivity Problems with Multiple particles, preprint available at http://arxiv.org/PS_cache/arxiv/pdf/0909/0909.3901v1.pdf

[7] Bonnetier, E., Vogelius, M. : An Elliptic Regularity Result for a Composite Medium with "Touching" Fibers of Circular Cross-Section, *SIAM J. Math. Anal.*, **31:3**, 2000, pp. 651–677.

[8] Einziner R.: Metal oxide varistors, *Annual Review of Materials Research*, **17**, 1987, pp. 299–321.

[9] Evans, L. C. : Partial Differential Equations, *American Mathematical Society, Graduate Studies in Mathematics*, **19**, 1998.

[10] Evans, L. C. amd Gangbo, W. : Differential equations methods for the Monge-Kantorovich mass transfer problem, *Mem. Amer. Math. Soc.*, **137:653**, 1999.

[11] Garroni, A., Nesi, V., Ponsiglione, M. : Dielectric breakdown: optimal bounds. *R. Soc. Lond. Proc. Ser. A Math. Phys. Eng. Sci.*, **457(2014)**, 2001, pp. 2317-2335.

[12] Huang, J. P., Yu, K. W. : Effective Nonlinear Optical Properties of Graded MetalDielectric Composite Films of Anisotropic Particles, *J. Opt. Soc. Am. B*, **22:8**, 2005, pp. 1640–1647.

[13] Idiart, M. : The macroscopic behavior of power-law and ideally plastic materials with elliptical distribution of porosity, *Mechanics Research Communications*, **35**, 2008, pp.583–588.

[14] Keller, J.B.: Conductivity of a Medium Containing a Dense Array of Perfectly Conducting Spheres or Cylinders or Nonconducting Cylinders, *J. Appl. Phys.*, **34:4**, 1963, pp. 991–993.

[15] Levy, O., Kohn, R. V. : Duality relations for non-Ohmic composites, with applications to behavior near percolation. *J. Statist. Phys.*, **90:1-2**, 1998, pp. 159-189.

[16] Ladyzhenskaya, O. A., Uraltseva, N. N. : *Linear and quasilinear equations of elliptic type*. Second edition, revised. Izdat. "Nauka", Moscow, 1973.

[17] Li, Y. Y., Nirenberg, L. : Estimates for Ellliptic System from Composite Material, *Comm. Pure Appl. Math.*, **56:7**, 2003, pp. 892–925.

[18] Li, Y. Y., Vogelius, M. : Gradient Estimates for Solution to Divergence Form Elliptic Equation with Discontinuous Coefficients, *Arch. Rational Mech. Anal.*, **153**, 2000, pp. 91–151.

[19] Lieberman, G. M. : Boundary regularity for solutions of degenerate elliptic equations, *Nonlinear Anal.*, **12:11**, 1988, pp. 1203–1219.

[20] Lim, M., Yun, K. : Blow-up of Electric Fields between Closely Spaced Spherical Perfect Conductors, *Comm. in PDEs*, **34:10**, 2009, pp. 1287–1315.

[21] Lipton, R.: Optimal Lower Bounds on the Electric Field Concentration in Composite Media, *Journal of Applied Physics*, **96**, 2004, pp. 2821–2827.

[22] Novikov A. : A Discrete Network Approximation for Effective Conductivity of non-Ohmic High Contast Composites, *Commun. Math. Sci.*, **7:3**, 2009, pp. 719–740.

[23] Ponte Castaneda, P., Suquet, P. : Nonlinear composties. *Advances in Applied Mechanics*, **34**, 1997, pp. 171-302.

[24] Ružička, M. : *Electrorheological fluids: modeling and mathematical theory*, Lecture Notes in Mathematics, **1748**. Springer-Verlag, Berlin, 2000.

[25] Le, V. K., Schmitt, K. : Sub-Supersolution Theorems for Quasilinear Elliptic Problems: A Variational Approach, *Electronic Journal of Differential Equations*, **118**, 2004, pp. 1–7.

[26] Suquet, P. : Overall potentials and extremal surfaces of power law or ideally plastic composites, *J. Mech. Phys. Solids*, **41(6)**, 1993, pp.981-1002.

[27] Talbot, D. R. S., Willis, J. R. : Upper and lower bounds for the overall properties of a nonlinear elastic composite dielectric. I. Random microgeometry. *Proc. R. Soc. Lond., A*, **447**, 1994, pp.365-396.

[28] Tolksdorf, P. : Regularity for a More General Class of Quasilinear Elliptic Equations, *J. Differential Equations*, **51**, 1984, pp.126–150.

[29] Yun, K. : Estimates for Electric Fields Blown Up Between Closely Adjacent Conductors with Arbitrary Shape, *SIAM J. Appl. Math.*, **67:3**, 2007, pp. 714–730.

[30] Yun, K. : Optimal bound on high stresses occurring between stiff fibers with arbitrary shaped cross-sections, *J. Math. Anal. Appl.*, **350**, 2009, pp. 306–312.